\documentclass{article}
\usepackage{amsmath,amsthm,amsopn,amstext,amscd,amsfonts,amssymb}
\usepackage{natbib}

\def\tr{\mathop{\rm tr}\nolimits}
\def\rank{\mathop{\rm rank}\nolimits}
\def\build#1#2#3{\mathrel{\mathop{#1}\limits^{#2}_{#3}}}

\def\vech{\mathop{\rm vech}\nolimits}

\def\etr{\mathop{\rm etr}\nolimits}

\newcommand {\findemo}{\hfill \, $\Box$ \,\\[2ex]}

\renewenvironment{abstract}
                 {\vspace{6pt}
                  \begin{center}
                  \begin{minipage}{5in}
                  \centerline{\textbf{Abstract}}
                  \noindent\ignorespaces
                 }
                 {\end{minipage}\end{center}}
\newtheorem{thm}{\textbf{Theorem}}[section]
\newtheorem{cor}{\textbf{Corollary}}[section]
\newtheorem{lem}{\textbf{Lemma}}[section]
\theoremstyle{definition}

\title{\Large \textbf{Shape Theory via QR decomposition}}
\author{
  \textbf{Jos\'e A. D\'{\i}az-Garc\'{\i}a} \thanks{Corresponding author\newline
    }\\
  {\normalsize Department of Statistics and Computation} \\
  {\normalsize Universidad Aut\'onoma Agraria Antonio Narro}\\
  {\normalsize 25350 Buenavista, Saltillo, Coahuila, Mexico} \\
  {\normalsize E-mail: jadiaz@uaaan.mx} \\[2ex]
  \textbf{Francisco J. Caro-Lopera} \\
  {\normalsize Department of Basic Sciences} \\
  {\normalsize Universidad de Medell\'{\i}n} \\
  {\normalsize Carrera 87 No.30-65, of. 5-103}\\
  {\normalsize Medell\'{\i}n, Colombia}\\
  {\normalsize E-mail: fjcaro@udem.edu.co}\\
}
\date{}
\begin{document}
\maketitle

\begin{abstract}
This work sets the non isotropic noncentral elliptical shape distributions via QR
decomposition in the context of zonal polynomials, avoiding the invariant polynomials
and the open problems for their computation. The new shape distributions are easily
computable and then the inference procedure can be studied under exact densities
instead under the published approximations and asymptotic densities under isotropic
models. An application in Biology is studied under the classical gaussian approach
and a two non gaussian models.
\end{abstract}

\section{Introduction and the main principle}\label{Sec:Introd}

Considering the non isotropy in the non central shape theory has been very
problematic, even in the gaussian case, see for example \citet{GM93}; the
corresponding shape densities involve expansions products of powers of traces of
different matrices, which forces the apparition of invariant polynomials
(\citet{d:80}). Then the resulting densities enlarge the list of uncountable
densities in the noncentral multivariate statistics, which can not be computable, and
remains as theoretical results, very far from the inference and the applications. So,
the applications in shape theory have been force to avoid those polynomials, but at a
very high cost, the assumption of isotropy. However, in this, the resulting densities
were series of zonal polynomials, and they could not studied properly (before the
works of \citet{KE06}), and again they forced to the use of approximations and
asymptotic distribution to perform inference, \citet{GM93}, \citet{DM98}, and the
references therein.

The following principle solves the first and more important problem, avoiding the
invariant polynomials, and setting the new shape distributions in terms of series of
zonal polynomials, this series can be computed accurately and efficiently by simple
modifications of the powerful algorithms of hypergeometric functions given by
\citet{KE06}.

From the point of view of applications, the isotropic assumption
$\mathbf{\Theta}=\mathbf{I}_{K}$ for an elliptical shape model of the form
$\mathbf{X} \sim \mathcal{E}_{N \times K} (\boldsymbol{\mu}_{{}_{\mathbf{X}}},
\mathbf{\Sigma}_{{}_{\mathbf{X}}}, \mathbf{\Theta},h)$, restricts substantially the
correlations of the landmarks in the figure.  So, we expect the non isotropic model,
with any positive definite matrix $\mathbf{\Theta}$, as the best model for
considering all the possible correlations among the anatomical (geometrical o
mathematical) points. However, using the classical approach of the published
literature of shape (see for example \citet{GM93}) under the non isotropic model, we
obtain immediately invariant polynomials, which can not be computed at this time for
large degrees.

In order to avoid this problem consider the following procedure: Let
$$
  \mathbf{X} \sim \mathcal{E}_{N \times K} (\boldsymbol{\mu}_{{}_{\mathbf{X}}},
  \mathbf{\Sigma}_{{}_{\mathbf{X}}}, \mathbf{\Theta},h),
$$
if $\mathbf{\Theta}^{1/2}$ is the positive definite square root of the matrix
$\mathbf{\Theta}$, i .e. $\mathbf{\Theta} = (\mathbf{\Theta}^{1/2})^{2}$, with
$\mathbf{\Theta}^{1/2}:$ $K \times K$, \citet[p. 11]{gv:93}, and noting that
$$
  \mathbf{X} \mathbf{\Theta}^{-1} \mathbf{X}' = \mathbf{X} (\mathbf{\Theta}^{-1/2}
  \mathbf{\Theta}^{-1/2})^{-1}\mathbf{X}' = \mathbf{X} \mathbf{\Theta}^{-1/2}
  (\mathbf{X} \mathbf{\Theta}^{-1/2})' = \mathbf{Z}\mathbf{Z}',
$$ where
$$
\mathbf{Z} = \mathbf{X} \mathbf{\Theta}^{-1/2},
$$
then
$$\mathbf{Z} \sim \mathcal{E}_{N \times K}(\boldsymbol{\mu}_{{}_{\mathbf{Z}}}, \mathbf{\Sigma}_{{}_{\mathbf{X}}}, \mathbf{I}_{K},
h)$$
 with $\boldsymbol{\mu}_{{}_{\mathbf{Z}}} = \boldsymbol{\mu}_{{}_{\mathbf{X}}} \mathbf{\Theta}^{-1/2}$, (see \citet[p. 20]{gv:93}).

And we arrive at the classical starting point in shape theory where the original
landmark matrix is replaced by $\mathbf{Z} = \mathbf{X} \mathbf{\Theta}^{-1/2}$. Then
we can proceed as usual, removing from $\mathbf{Z}$, translation, scale, rotation
and/or reflection in order to obtain the shape of $\mathbf{Z}$ (or $\mathbf{X}$) via
the QR decomposition, for example.

Namely, the QR shape coordinates $\mathbf{u}$ of $\mathbf{X}$ are constructed in
several steps summarized in the expression
\begin{equation}\label{eq:QRSteps}
    \mathbf{L}\mathbf{X}\mathbf{\Theta}^{-1/2}=\mathbf{L}\mathbf{Z}=\mathbf{Y}=\mathbf{TH}=
    r\mathbf{WH}=r\mathbf{W}(\mathbf{u})\mathbf{H},
\end{equation}
which we discuss next. Observe that $\boldsymbol{\mu}_{{}_{\mathbf{Z}}} =
\boldsymbol{\mu}_{{}_{\mathbf{X}}} \mathbf{\Theta}^{-1/2}$ and the QR shape
coordinates of $\boldsymbol{\mu}_{{}_{\mathbf{Z}}}$ are defined analogously. The
matrix $\mathbf{L}:(N-1)\times N$ has orthonormal rows to $\textbf{1}=(1,\ldots,1)'$.
$\mathbf{L}$ can be a submatrix of the Helmert matrix, for example.

Let $\boldsymbol{\mu}=\mathbf{L}\boldsymbol{\mu}_{{}_{\mathbf{X}}}$, then
$\mathbf{Y}:(N-1)\times K$ is invariant to translations of the figure $\mathbf{Z}$,
and
$$
  \mathbf{Y}\sim\mathcal{E}_{N-1 \times K}(\boldsymbol{\mu}\mathbf{\Theta}^{-1/2},
  \mathbf{\Sigma}\otimes \mathbf{I}_{K},h),
$$
where $\mathbf{\Sigma}=\mathbf{L}\mathbf{\Sigma}_{\mathbf{X}}\mathbf{L}'$.

Now, let be $n=\min(N-1,K)$ and $p=\rank \mu$. In (\ref{eq:QRSteps}), $\mathbf{Y} =
\mathbf{TH}$ is the QR decomposition, where $\mathbf{T}:(N-1)\times n$ is lower
triangular with $t_{ii}>0$, $i=1,\ldots, \min(n,K-1)$, and $\mathbf{H}:n\times K$,
$\mathbf{H}\in \mathcal{V}_{n,K}$, the Stiefel manifold. Note that $\mathbf{T}$ is
invariant to translations and rotations of $\mathbf{Z}$. The matrix $\mathbf{T}$ is
referred as the \textit{QR size-and-shape} and their elements are the QR
size-and-shape coordinates of the original landmark data $\mathbf{Z}$. Typically in
shape analysis there are more landmarks than dimensions ($N>K$). $\mathbf{H}$ acts on
the right to transform $\Re^{K}$ instead of acting on the left as in the multivariate
analysis. In our case we see the landmarks as variables and the dimensions as
observations, then the transposes of our matrices $\mathbf{Z}$ and $\mathbf{Y}$ can
be seen as classical multivariate data matrices.

According to the nature of the base $\mathbf{H}$ and providing that $N-1\geq K$, we
say that $\mathbf{T}$ contains the \textit{QR reflection size-and-shape} coordinates
if $\mathbf{H}$ includes reflection, i.e. $\mathbf{H}\in \mathcal{O}(K)$,
$|\mathbf{H}|=\pm 1 $ and $t_{KK}\geq 0$; otherwise, if $\mathbf{H}$ excludes
reflection, $\mathbf{H}\in \mathcal{SO}(K)$, $|\mathbf{H}|=+1$, $t_{KK}$ is not
restricted, we say that $\mathbf{T}$ contains the QR size-and-shape coordinates.
These cases will denote by $\mathbf{T}^{R}$ and $\mathbf{T}^{NR}$, respectively. In
the classical multivariate case, $n<K$, we do not have such classifications for
$\mathbf{T}$.

Now, if we divide $\mathbf{T}$ by its size, the centroid size of $\mathbf{Z}$,
$$
r=\|\mathbf{T}\|=\sqrt{\tr \mathbf{T}'\mathbf{T}}=\|\mathbf{Y}\|.
$$
we obtain the so called \textit{QR shape} matrix $\mathbf{W}$ in (\ref{eq:QRSteps}).
We define $\mathbf{W}^{R} = \mathbf{T}^{R}/r$ or $\mathbf{W}^{NR}=\mathbf{T}^{NR}/r$
if W includes or excludes reflection, respectively, and given that
$\|\mathbf{W}\|=1$, the elements of $\mathbf{W}$ are a direction vector for shape,
and $\mathbf{u}$ comprises $m=(N-1)K-nK+\frac{1}{2}n(n+1)-1$ generalized polar
coordinates.

Before deriving the main results of this paper we must solve some discrepancies
between the shape theory and the classical multivariate theory. Recall that for a
given $\mathbf{Y}:n\times K$, $n=N-1\leq K$, then $\mathbf{YY}'$ has the noncentral
Wishart distribution which is invariant to orientation and reflection, but if $n\geq
K$ that density does not exist with respect to the Lebesgue measure defined on the
space of positive definite $n \times n$ matrices, and we therefore use the
size-and-shape matrix $\mathbf{T}$. However, the density of $\mathbf{YY}'$ when,
$n\geq K$, exist on the $(nK-K(K-1)/2)$-dimensional manifold of rank-K positive
semidefinite $n \times n$ matrices with $K$ distinct positive eigenvalues, see
\citet{dggg05} and \citet{dggj06}. This last fact can provide an alternative form to
study the shape theory, which is being analysed by the authors at present.

And finally, classical integration over $\mathcal{O}(K)$ involving zonal polynomials
gives the density of $\mathbf{T}^{R}$, but $\mathbf{T}^{NR}$ demands integration over
$\mathcal{SO}(K)$, then we just recall that the corresponding integrals are the same
when $n<K$, and, for $n\geq K$ and $p<K$, they are twice the integral over
$\mathcal{SO}(K)$.

This work is distributed as follows: first, the size and shape distribution for any
elliptical model with a full Kronecker covariance matrix is derived in section
\ref{sec:QRsizeandshape}. The the shape density is obtained in section
\ref{sec:QRshape} and the classical isotropic gaussian shape density, full derived in
\citet{GM93}, follows here as a corollary, then the section
\ref{sec:QRexcludingreflections} describes the excluding reflection shape densities.
The central case of the shape density is studied in section \ref{sec:QRcentralcase},
and a remarkable property is established, i.e.  it is established that the central QR
reflection shape density is invariant under the elliptical family.  Finally, some
particular elliptical densities are derived in section \ref{sec:particularmodels} in
order to perform inference on exact distributions; a subfamily of Kotz distributions
which contains the gaussian one is derived, then two elements of that class (the
gaussian and a non gaussian model) is applied to an existing publish data, the mouse
vertebra study. Some test for detecting shape differences are gotten and the models
are discriminated by the use of a dimension criterion such as the modified BIC
criterion.

\section{QR Size-and-shape distribution}\label{sec:QRsizeandshape}

\begin{lem}\label{lem:jacobianQR}
Let $\mathbf{Y}:(N-1)\times K$, then there exists a $\mathbf{T}:(N-1)\times n$ lower
triangular matrix with $t_{ii}\geq 0$, $i=1,\ldots ,\min (n,K-1)$, and $\mathbf{H}\in
\mathcal{V}_{n,k}$ such that $\mathbf{Y}=\mathbf{T}\mathbf{H}$ and
\begin{equation}\label{eq:jacobianQR}
(d\mathbf{Y})=\prod_{i=1}^{n}t_{ii}^{K-i}(d\mathbf{T})(\mathbf{H}d\mathbf{H}')
\end{equation}
\end{lem}

\begin{lem}\label{lem:James64eq22}
Let $\mathbf{A}:r\times s$, and $\mathbf{H}\in \mathcal{V}_{s,m}$ then
\begin{equation}\label{eq:James64eq22}
\int_{\mathbf{H}\in V_{s,m}}(\tr
\mathbf{A}\mathbf{H})^{2t}(\mathbf{H}d\mathbf{H}')=\frac{2^{s}\pi^{sm/2}}{\Gamma_{s}
\left[\frac{1}{2}m\right]}\sum_{\kappa}\frac{\left(\frac{1}{2}\right)_{t}}
{\left(\frac{1}{2}m\right)_{\kappa}}C_{\kappa}(\mathbf{AA}'),
\end{equation}
\end{lem}
where $C_{\kappa}(\mathbf{B})$ are the zonal polynomials of $\mathbf{B}$
corresponding to the partition $\kappa=(t_{1},\ldots t_{\alpha})$ of $t$, with
$\sum_{i=1}^{\alpha}t_{i}=t$; and $(a)_{\kappa}=\prod_{i=1}(a-(j-1)/2)_{t_{j}}$,
$(a)_{t}=a(a+1)\cdots (a+t-1)$, are the generalized hypergeometric coefficients and
$\Gamma_{s}(a)=\pi^{s(s-1)/4}\prod_{j=1}^{s}\Gamma(a-(j-1)/2)$ is the multivariate
Gamma function.

\textit{Proof.}  It is follows from \citet{JAT64}, eq. (22) and \citet{MR1982}, lemma
9.5.3, p. 397.\findemo

\begin{thm}\label{th:QRreflectionsizeandshape}
The QR reflection size-and-shape is
\begin{equation}\label{eq:QRreflectionsizeandshape}
f_{\mathbf{T}}(\mathbf{T})=\frac{2^{n}\pi^{nK/2}\displaystyle\prod_{i=1}^{n}t_{ii}^{K-i}}
{\Gamma_{n}\left[\frac{1}{2}K\right]|\mathbf{\Sigma}|^{K/2}}
\sum_{t=0}^{\infty}\sum_{\kappa}\frac{h^{(2t)}[\tr
(\mathbf{\Sigma}^{-1}\mathbf{T}\mathbf{T}'+\mathbf{\Omega})]}{t!}\frac{C_{\kappa}(\mathbf{\Omega}
\mathbf{\Sigma}^{-1}\mathbf{T}\mathbf{T}')}{\left(\frac{1}{2}K\right)_{\kappa}},
\end{equation}
\end{thm}
where $\mathbf{\Omega} = \mathbf{\Sigma}^{-1} \boldsymbol{\mu} \mathbf{\Theta}^{-1}
\boldsymbol{\mu}'$, $C_{\kappa}(\mathbf{B})$ are the zonal polynomials of
$\mathbf{B}$ corresponding to the partition $\kappa=(t_{1},\ldots, t_{\alpha})$ of
$t$, with $\sum_{i=1}^{\alpha}t_{i}=t$ and $h^{(j)}(v)$ is the $j$-th derivative of
$h$ with respect to $v$.

\textit{Proof.}  The density of $\mathbf{Y}$ is
\begin{eqnarray*}
f_{\mathbf{Y}}(\mathbf{Y})&=&\frac{1}{|\mathbf{\Sigma}|^{K/2}}h\left[\tr \mathbf{\Sigma}^{-1}
(\mathbf{Y}-\boldsymbol{\mu}\mathbf{\Theta}^{-1/2})(\mathbf{Y}-\boldsymbol{\mu}\mathbf{\Theta}^{-1/2})'\right]\\
&=&\frac{1}{|\mathbf{\Sigma}|^{K/2}}h\left[\tr
\left(\mathbf{\Sigma}^{-1}\mathbf{Y}\mathbf{Y}'
+\mathbf{\Sigma}^{-1} \boldsymbol{\mu} \mathbf{\Theta}^{-1} \boldsymbol{\mu}' - 2
\mathbf{\Sigma}^{-1}\mathbf{Y}\mathbf{\Theta}^{-1/2}\boldsymbol{\mu}'\right)\right]\\
&=&\frac{1}{|\mathbf{\Sigma}|^{K/2}}h\left[\tr
\left(\mathbf{\Sigma}^{-1}\mathbf{Y}\mathbf{Y}'+\mathbf{\Omega}\right)-2\tr
\mathbf{\Sigma}^{-1}\mathbf{Y}\mathbf{\Theta}^{-1/2}\boldsymbol{\mu}'\right],
\end{eqnarray*}
with $\mathbf{\Omega} = \mathbf{\Sigma}^{-1} \boldsymbol{\mu} \mathbf{\Theta}^{-1}
\boldsymbol{\mu}' $. Taking $\mathbf{Y}=\mathbf{T}\mathbf{H}$ and using Lemma
\ref{lem:jacobianQR}, the joint density of $\mathbf{\mathbf{H}}$ and $\mathbf{T}$ is
$$
  f_{\mathbf{H},\mathbf{T}}(\mathbf{H},\mathbf{T})=\frac{\displaystyle\prod_{i=1}^{n}t_{ii}^{K-i}}
  {|\mathbf{\Sigma}|^{K/2}}h\left[\tr \left(\mathbf{\Sigma}^{-1}\mathbf{T}\mathbf{T}'+\mathbf{\Omega}\right)-2\tr
  \mathbf{\Theta}^{-1/2}\boldsymbol{\mu}'\mathbf{\Sigma}^{-1}\mathbf{T}\mathbf{H}\right].
$$
Assuming that $h(\cdot)$ can be expanded in power series, see
\citet{fz:90}, i.e.
$$
h(a+v)=\sum_{t=0}^{\infty}\frac{h^{t}(a)v^{t}}{t!}.
$$
Thus
$$
  f_{\mathbf{H},\mathbf{T}}(\mathbf{H},\mathbf{T})=\frac{\displaystyle\prod_{i=1}^{n}t_{ii}^{K-i}}
  {|\mathbf{\Sigma}|^{K/2}}\sum_{t=0}^{\infty}\frac{1}{t!}h^{t}\left[\tr
  \left(\mathbf{\Sigma}^{-1}\mathbf{T}\mathbf{T}'+\mathbf{\Omega}\right)\right]\left[\tr
  \left(-2\mathbf{\Theta}^{-1/2}\boldsymbol{\mu}'\mathbf{\Sigma}^{-1}\mathbf{T}\mathbf{H}\right)\right]^{t}.
$$
Now, for integration on $\mathbf{H}\in \mathcal{V}_{n,K}$, we note that it is zero
when $t$ is odd (see Theorem II, p.876 in James (1961) or eqs. (44)-(46) in James
(1964)).  The marginal density of $\mathbf{T}$ is expressed as
\begin{eqnarray*}
f_{\mathbf{T}}(\mathbf{T})&=&\frac{\displaystyle\prod_{i=1}^{n}t_{ii}^{K-i}}{|\mathbf{\Sigma}|^{K/2}}
\sum_{t=0}^{\infty}\frac{1}{(2t)!}h^{2t}\left[\tr \mathbf{\Sigma}^{-1}
\left(\mathbf{T}\mathbf{T}'+\mathbf{\Omega}\right)\right]\\
&&\times \int_{\mathcal{V}_{n,K}}\left[\tr
\left(-2\mathbf{\Theta}^{-1/2}\boldsymbol{\mu}'
\mathbf{\Sigma}^{-1}\mathbf{T}\mathbf{H}\right)\right]^{2t}(\mathbf{H}d\mathbf{H}').
\end{eqnarray*}
So, by Lemma \ref{lem:James64eq22} and recalling that
$C_{\lambda}(a\mathbf{A})=a^{t}C_{\lambda}(\mathbf{A})$, for $a$ constant, we have
$$
  \int_{\mathcal{V}_{n,K}}\left[\tr \left(-2\mathbf{\Theta}^{-1/2}\boldsymbol{\mu}'
\mathbf{\Sigma}^{-1}\mathbf{T}\mathbf{H}\right)\right]^{2t}(\mathbf{H}d\mathbf{H}')
\hspace{6cm}
$$

\vspace{-.5cm}

\begin{eqnarray*}
\hspace{3cm}&=&\frac{2^{n}\pi^{nK/2}}{\Gamma_{n}\left[\frac{1}{2}K\right]}
\sum_{\kappa}\frac{\left(\frac{1}{2}\right)_{t}}{\left(\frac{1}{2}K\right)_{\kappa}}
C_{\kappa}\left(4\mathbf{\Theta}^{-1/2}\boldsymbol{\mu}'\mathbf{\Sigma}^{-1}\mathbf{T}\mathbf{T}'
\mathbf{\Sigma}^{-1}\boldsymbol{\mu}\mathbf{\Theta}^{-1/2}\right)\\
&=&\frac{2^{n}\pi^{nK/2}}{\Gamma_{n}\left[\frac{1}{2}K\right]}
\sum_{\kappa}\frac{\left(\frac{1}{2}\right)_{t}4^{t}}{\left(\frac{1}{2}K\right)_{\kappa}}
C_{\kappa}\left(\mathbf{\Theta}^{-1}\boldsymbol{\mu}'\mathbf{\Sigma}^{-1}\mathbf{T}\mathbf{T}'
\mathbf{\Sigma}^{-1}\boldsymbol{\mu}\right)\\
&=&\frac{2^{n}\pi^{nK/2}}{\Gamma_{n}\left[\frac{1}{2}K\right]}
\sum_{\kappa}\frac{\left(\frac{1}{2}\right)_{t}4^{t}}{\left(\frac{1}{2}K\right)_{\kappa}}
C_{\kappa}\left(\mathbf{\Omega}\mathbf{\Sigma}^{-1}\mathbf{T}\mathbf{T}'\right).
\end{eqnarray*}
From \citet{MR1982}, p.21
$\frac{\Gamma\left(k+\frac{1}{2}n\right)}{\Gamma\left(\frac{1}{2}n\right)}=\left(\frac{1}{2}n\right)_{k}$
then
$\frac{\Gamma\left(k+\frac{1}{2}\right)}{\Gamma\left(\frac{1}{2}\right)}=\left(\frac{1}{2}\right)_{k}$
and $\frac{\left(\frac{1}{2}\right)_{k}4^{k}}{(2k)!}=\frac{1}{k!}$, in our case
$k=t$, and the result follows.  \findemo

Alternatively, the size-and-shape density (\ref{eq:QRreflectionsizeandshape}) can be
obtained as a particular case of the singular case studied in \citet{dggg05}.

\section{QR Shape distribution}\label{sec:QRshape}
Now, observe that for $\mathbf{T}:(N-1)\times n$, $n=\min (N-1,K)$, the matrix
$\mathbf{T}$ contains $(N-1)K-nK+n(n+1)/2$ non null QR rectangular coordinates
$(t_{ij}\neq 0)$. Let $\vech \mathbf{T}$ a vector consisting of the no null elements
of $\mathbf{T}$, taken column by column. Then the QR shape matrix $\mathbf{W}$ can be
written as
$$
\vech \mathbf{W}=\frac{1}{r}\vech \mathbf{T},\quad r=||\mathbf{T}||=\sqrt{\tr
\mathbf{T}'\mathbf{T}}=||\mathbf{Y}||,
$$
then by Theorem 2.1.3, p.55 of \citet{MR1982},
$$
  (d \vech \mathbf{T}) = r^{m}\prod_{i=1}^{m}\sin^{m-i}\theta_{i}
  \left(\bigwedge_{i=1}^{m}d\theta_{i}\right)\wedge
dr,
$$
with $m=(N-1)K-nK+n(n+1)/2-1$.  Denoting $\mathbf{u} = (\theta_{1}, \ldots,
\theta_{m})'$ and $J(\mathbf{u}) = r^{m}\prod_{i=1}^{m}\sin^{m-i}\theta_{i}$, so
$$
(d\mathbf{T})=r^{m}J(\mathbf{u})\left(\bigwedge_{i=1}^{m}d\theta_{i}\right)\wedge dr.
$$

\begin{thm}\label{th:QRreflectionshape}
The QR reflection shape density is
\begin{eqnarray}\label{eq:QRreflectionshape}
f_{\mathbf{W}}(\mathbf{W})&=&\frac{2^{n}\pi^{nK/2}\displaystyle
\prod_{i=1}^{n}w_{ii}^{K-i}J(\mathbf{u})}
{\Gamma_{n}\left[\frac{1}{2}K\right]|\mathbf{\Sigma}|^{K/2}}
\sum_{t=0}^{\infty}\sum_{\kappa}\frac{C_{\kappa}(\mathbf{\Omega}
\mathbf{\Sigma}^{-1}\mathbf{W}\mathbf{W}')}{t!\left(\frac{1}{2}K\right)_{\kappa}}\nonumber\\
&&\times\int_{0}^{\infty} r^{M+2t-1}h^{(2t)}[r^{2}\tr
\mathbf{\Sigma}^{-1}\mathbf{W}\mathbf{W}'+\tr\mathbf{\Omega}](dr),
\end{eqnarray}
\end{thm}
where $M=(N-1)K$.

\textit{Proof.}  The density of $\mathbf{T}$ is
\begin{equation*}
f_{\mathbf{T}}(\mathbf{T})=\frac{2^{n}\pi^{nK/2}}{\Gamma_{n}[K/2]}\frac{\displaystyle
\prod_{i=1}^{n}t_{ii}^{K-i}}
{|\mathbf{\Sigma}|^{K/2}}\sum_{t=0}^{\infty}\sum_{\kappa}
\frac{h^{(2t)}[\tr(\mathbf{\Sigma}^{-1}\mathbf{T}\mathbf{T}'+\mathbf{\Omega})]}{t!}
\frac{C_{\kappa}(\mathbf{\Omega}\mathbf{\Sigma}^{-1}\mathbf{T}\mathbf{T}')}{\left(\frac{1}{2}K\right)_{\kappa}}.
\end{equation*}
Making the change of variables $\mathbf{W}(\mathbf{u})=\mathbf{T}/r$, the joint
density function of $r$ and $\mathbf{u}$ is
\begin{eqnarray*}
f_{r,\mathbf{W}}(r,\mathbf{W})&=&\frac{2^{n}\pi^{nK/2}}{\Gamma_{n}[K/2]}\frac{\displaystyle
\prod_{i=1}^{n}
(rw_{ii})^{K-i}}{|\mathbf{\Sigma}|^{K/2}}\sum_{t=0}^{\infty}\sum_{\kappa}
\frac{h^{(2t)}[\tr(r^{2}\mathbf{\Sigma}^{-1}\mathbf{W}\mathbf{W}'+\mathbf{\Omega})]}{t!}
\\&&\times\frac{C_{\kappa}(r^{2}\mathbf{\Omega}\mathbf{\Sigma}^{-1}\mathbf{W}\mathbf{W}')}
{\left(\frac{1}{2}K\right)_{\kappa}}r^{m}J(\mathbf{u}).
\end{eqnarray*}
Now, note that
\begin{itemize}
    \item
    $C_{\kappa}(r^{2}\mathbf{\Omega}\mathbf{\Sigma}^{-1}\mathbf{W}\mathbf{W}')=r^{2t}
    C_{\kappa}(\mathbf{\Omega}\mathbf{\Sigma}^{-1}\mathbf{W}\mathbf{W}')$.
    \item $\displaystyle
    \displaystyle \prod_{i=1}^{n}(rw_{ii})^{K-i}=\sum_{i=1}^{n}r^{(K-i)}\displaystyle \prod_{i=1}^{n}w_{ii}^{K-i}=r^{nK -
    \frac{n(n+1)}{2}}\displaystyle \prod_{i=1}^{n}w_{ii}^{K-i}.
    $
    \item
    $h^{(2t)}[\tr(r^{2}\mathbf{\Sigma}^{-1}\mathbf{W}\mathbf{W}'+\mathbf{\Omega})] =
    h^{(2t)}[r^{2}\tr\mathbf{\Sigma}^{-1}\mathbf{W}\mathbf{W}'+\tr\mathbf{\Omega}]$.
\end{itemize}
Collecting powers of $r$ by $r^{nK-\frac{n(n+1)}{2}+2t+m}=r^{M+2t-1}$, $M=(N-1)K$,
the marginal of $\mathbf{W}$ is:
\begin{eqnarray*}
f_{\mathbf{W}}(\mathbf{W})&=&\frac{2^{n}\pi^{nK/2}\displaystyle
\prod_{i=1}^{n}w_{ii}^{K-i}J(\mathbf{u})}
{\Gamma_{n}\left[\frac{1}{2}K\right]|\mathbf{\Sigma}|^{K/2}}
\sum_{t=0}^{\infty}\sum_{\kappa}\frac{C_{\kappa}(\mathbf{\Omega}
\mathbf{\Sigma}^{-1}\mathbf{W}\mathbf{W}')}{t!\left(\frac{1}{2}K\right)_{\kappa}}\nonumber\\
&&\times\int_{0}^{\infty}r^{M+2t-1}h^{(2t)}[r^{2}\tr
\mathbf{\Sigma}^{-1}\mathbf{W}\mathbf{W}' + \tr\mathbf{\Omega}](dr). \qed
\end{eqnarray*}

When $\mathbf{\Sigma} = \sigma^{2}\mathbf{I}$, then $\mathbf{\Omega} =
\boldsymbol{\mu} \mathbf{\Theta}^{-1} \boldsymbol{\mu}'/\sigma^{2}$,
$|\mathbf{\Sigma}|^{K/2} = \sigma^{M}$, and $r^{2}\tr
\mathbf{\Sigma}^{-1}\mathbf{W}\mathbf{W}' = r^{2}/\sigma^{2}$, because $\tr
\mathbf{W}\mathbf{W}' = 1$, thus Theorem \ref{th:QRreflectionshape} becomes

\begin{cor}\label{coroSHdiagonal}
The QR reflection shape density is
\begin{eqnarray}\label{eq:QRreflectionshapediagonal}
f_{\mathbf{W}}(\mathbf{W})&=&\frac{2^{n}\pi^{nK/2}\displaystyle
\prod_{i=1}^{n}w_{ii}^{K-i}J(\mathbf{u})}
{\Gamma_{n}\left[\frac{1}{2}K\right]\sigma^{M}}
\sum_{t=0}^{\infty}\sum_{\kappa}\frac{C_{\kappa}\left(\frac{1}{\sigma^{2}}\mathbf{\Omega}
\mathbf{W}\mathbf{W}'\right)}{t!\left(\frac{1}{2}K\right)_{\kappa}}\nonumber\\
&& \times \int_{0}^{\infty}
r^{M+2t-1}h^{(2t)}[r^{2}/\sigma^{2}+\tr\mathbf{\Omega}](dr);
\end{eqnarray}
\end{cor}
and, if the gaussian model is considered with $\mathbf{\Theta}=\mathbf{I}$, the
resulting density corresponds with \citet[Theorem 2]{GM93}, see section
\ref{sec:particularmodels} below.

\section{Distributions excluding reflections}\label{sec:QRexcludingreflections}

From subsection 2.1 of \citet{GM93}, we can derive the QR size-and-shape and QR shape
densities excluding reflection:

\begin{itemize}
    \item If $n<K$, then Theorems \ref{th:QRreflectionsizeandshape},
\ref{th:QRreflectionshape} stand for the corresponding $\mathbf{T}=\mathbf{T}^{NR}$
and $\mathbf{W}=\mathbf{W}^{NR}$ excluding reflection densities.
    \item When $N-1\geq K$ and $p<K$, the QR size-and-shape density for
    $\mathbf{T}=\mathbf{T}^{NR}$  is
    (\ref{eq:QRreflectionsizeandshape}) divided by 2, where $t_{ii}\geq
    0$, for $i=1,\ldots,K-1$ and $t_{KK}$ is unrestricted.  When
    $N-1<K$ (\ref{eq:QRreflectionsizeandshape}) stands, since $t_{KK}$
    is not present.
    \item When $N-1\geq K$ and $p<K$, the QR shape density for $\mathbf{W}=\mathbf{W}^{NR}$
     is (\ref{eq:QRreflectionshape}) divided by 2,
    and $w_{KK}$ is unrestricted. When $N-1<K$
    (\ref{eq:QRreflectionshape}) holds, since $w_{KK}$ is not present.
    \item The preceding results also hold when $\rank \boldsymbol{\mu}=K$ and $\rank
    \mathbf{T}<K$, and event with probability zero.
    \item However, if $p=K$, the excluding reflection densities do
    not follow the above rule. For the  gaussian case,  see
    \citet{GM93} and \citet{gm:91}.
\end{itemize}

\section{Central Case}\label{sec:QRcentralcase}
The central case of the preceding sections can be derived easily.
\begin{cor}\label{cor:QRcentralreflectionsizeandshape}
The central QR reflection size-and-shape density is given by
\begin{equation*}
    f_{\mathbf{T}}(\mathbf{T})=\frac{2^{n}\pi^{\frac{nK}{2}}}{\Gamma_{n}\left[\frac{1}{2}K\right]
    |\mathbf{\Sigma}|^{\frac{K}{2}}}h[\tr\mathbf{\Sigma}^{-1}\mathbf{T}\mathbf{T}'].
\end{equation*}
\end{cor}
\textit{Proof.} It is straightforward from Theorem \ref{th:QRreflectionsizeandshape}
just take $\boldsymbol{\mu} = \mathbf{0}$ and recall that
$h^{(0)}[\tr\cdot]=h[\tr\cdot]$. \findemo

And:
\begin{cor}\label{cor:QRcentralreflectionshape}
The central QR reflection shape density is given by
\begin{equation*}
    f_{\mathbf{W}}(\mathbf{W})=\frac{2^{n}\pi^{\frac{nK}{2}}\displaystyle \prod_{i=1}^{n}w_{ii}^{K-i}J(\mathbf{u})}
    {\Gamma_{n}\left[\frac{1}{2}K\right]|\mathbf{\Sigma}|^{\frac{K}{2}}}
    \int_{0}^{\infty}r^{M-1}h[r^{2}\tr\mathbf{\Sigma}^{-1}\mathbf{W}\mathbf{W}'](dr).
\end{equation*}
\end{cor}
\textit{Proof.} Just take $\boldsymbol{\mu} = \mathbf{0}$ and
$h^{(0)}[\tr\cdot]=h[\tr\cdot]$ in Theorem \ref{th:QRreflectionshape}. \findemo

The corresponding central excluding reflection densities follows
according to Section \ref{sec:QRexcludingreflections}.

Observe that it is possible to obtain an invariant central shape density, i.e. the
density  does not depend on function $h(\cdot)$ Let $h$ be the density generator of
$\mathbf{Y}\sim \mathcal{E}_{N-1,K}(\mathbf{0},\mathbf{I} \otimes \mathbf{I},h)$,
i.e.
$$
f_{\mathbf{Y}}(\mathbf{Y})=h(\tr \mathbf{Y}\mathbf{Y}'),
$$
then by \citet{fz:90}, p.102, eq. 3.2.6, $$
\int_{0}^{\infty}r^{(N-1)K-1}h(r^{2})dr=\frac{\Gamma[(N-1)K/2]}{2\pi^{(N-1)K/2}}.
$$
So, if $s=(\tr\mathbf{\Sigma}^{-1}\mathbf{W}\mathbf{W}')^{1/2}r$,
$ds=(\tr\mathbf{\Sigma}^{-1}\mathbf{W}\mathbf{W}')^{1/2}(dr)$, then
\begin{eqnarray*}
&&\int_{0}^{\infty}r^{M-1}h[r^{2}\tr\mathbf{\Sigma}^{-1}\mathbf{W}\mathbf{W}'](dr)\\
&=&\int_{0}^{\infty}\left(\frac{s}{(\tr\mathbf{\Sigma}^{-1}\mathbf{W}\mathbf{W}')^{1/2}}\right)^{M-1}
h(s^{2})\frac{ds}{(\tr\mathbf{\Sigma}^{-1}\mathbf{W}\mathbf{W}')^{1/2}}\\
&=&(\tr\mathbf{\Sigma}^{-1}\mathbf{W}\mathbf{W}')^{-M/2}\int_{0}^{\infty}s^{M-1}h(s^{2})ds\\
&=&(\tr\mathbf{\Sigma}^{-1}\mathbf{W}\mathbf{W}')^{-M/2}\frac{\Gamma[M/2]}{2\pi^{M/2}}.
\end{eqnarray*}
Thus:
\begin{cor}\label{cor:QRcentralreflectionshapeinvariant}
When $\boldsymbol{\mu}=\mathbf{0}$ the QR reflection shape density is invariant under
the elliptical family and it is given by
\begin{equation*}
    f_{w}(w)=\frac{2^{n-1}\pi^{\frac{nK-M}{2}}\Gamma\left[M/2\right]}{\Gamma_{n}
    \left[\frac{1}{2}K\right]|\mathbf{\Sigma}|^{\frac{K}{2}}}
   \displaystyle \prod_{i=1}^{n}w_{ii}^{K-i}J(\mathbf{u})(\tr\mathbf{\Sigma}^{-1}\mathbf{W}\mathbf{W}')^{-\frac{M}{2}}.
\end{equation*}
\end{cor}

As in the noncentral case, if $\mathbf{\Sigma} =
\sigma^{2}\mathbf{I}$, then $|\mathbf{\Sigma}|^{\frac{K}{2}} =
\sigma^{M}$ and
$(\tr\mathbf{\Sigma}^{-1}\mathbf{W}\mathbf{W}')^{-\frac{M}{2}} =
\sigma^{M}$, thus:

\begin{cor}\label{cor:QRcentralreflectionshapeinvariantdiagonal}
When $\boldsymbol{\mu}=\mathbf{0}$ and $\mathbf{\Sigma} =
\sigma^{2}\mathbf{I}$ the QR reflection shape density is invariant under the
elliptical family and it is given by
\begin{equation*}
   f_{w}(w)=\frac{2^{n-1}\pi^{\frac{nK-M}{2}}\Gamma\left[M/2\right]}{\Gamma_{n}
   \left[\frac{1}{2}K\right]} J(\mathbf{u}) \displaystyle \prod_{i=1}^{n}w_{ii}^{K-i}.
\end{equation*}
\end{cor}

\section{Some particular models}\label{sec:particularmodels}

Finally, we give explicit shapes densities for some elliptical models.

The Kotz type I model is given by
$$
  h(y)=\frac{R^{\tau-1+\frac{K(N-1)}{2}}}{\Gamma\left(\frac{K(N-1)}{2}\right)}{\pi^{K(N-1)/2}
 \Gamma\left(\tau-1+\frac{K(N-1)}{2}\right)}y^{\tau-1}\exp(-Ry),
$$
So, the corresponding $k$-th derivative is
$$
\frac{d^{k}[y^{\tau-1}\exp\{-Ry\}]}{dy^{k}}=(-R)^{k}y^{\tau-1}\exp\{-Ry\}\left\{1+\sum_{m=1}^{k}\binom{k}{m}
\left[\prod_{i=0}^{m-1}(\tau-1-i)\right](-Ry)^{-m}\right\}.
$$

It is of interest the gaussian case, i.e. when $\tau=1$ and $R=\frac{1}{2}$, here the
derivation is straightforward from the general density.

The required derivative follows easily, it is,
$h^{(k)}(y)=\frac{R^{\frac{K(N-1)}{2}}}{\pi^{\frac{K(N-1)}{2}}}(-R)^{k}\exp(-Ry)$ and
\begin{eqnarray*}
   \int_{0}^{\infty} &&r^{M+2t-1}h^{(2t)}[r^{2}\tr \mathbf{\Sigma}^{-1}\mathbf{W}\mathbf{W}' +
   \tr\mathbf{\Omega}]dr\\&&=\frac{R^{t}}{2\pi^{\frac{M}{2}}}
   \exp(-R \tr\mathbf{\Omega})\left(\tr \mathbf{\Sigma}^{-1}\mathbf{W}\mathbf{W}'\right)^{-\frac{M}{2}-t}
   \Gamma\left(\frac{M}{2}+t\right).
\end{eqnarray*}

\begin{eqnarray*}
f_{\mathbf{W}}(\mathbf{W})&=&\frac{2^{n}\pi^{nK/2}J(\mathbf{u})\displaystyle
\prod_{i=1}^{n}w_{ii}^{K-i}}
{\Gamma_{n}\left[\frac{1}{2}K\right]|\mathbf{\Sigma}|^{K/2}}
\sum_{t=0}^{\infty}\sum_{\kappa}\frac{C_{\kappa}(\mathbf{\Omega}
\mathbf{\Sigma}^{-1}\mathbf{W}\mathbf{W}')}{t!\left(\frac{1}{2}K\right)_{\kappa}}\nonumber\\&&
\times\frac{R^{t}}{2\pi^{\frac{M}{2}}} \exp(-R \tr\mathbf{\Omega})\left(\tr
\mathbf{\Sigma}^{-1}\mathbf{W}\mathbf{W}'\right)^{-\frac{M}{2}-t}
\Gamma\left(\frac{M}{2}+t\right)\\
&=&\frac{\exp(-R \tr\mathbf{\Omega})\left(\tr
\mathbf{\Sigma}^{-1}\mathbf{W}\mathbf{W}'\right)^{-\frac{M}{2}}J(\mathbf{u})\displaystyle
\prod_{i=1}^{n}w_{ii}^{K-i}} {\pi^{\frac{M-nK}{2}}2^{-n+1}
\Gamma_{n}\left[\frac{1}{2}K\right]|\mathbf{\Sigma}|^{K/2}}
\\&&\times
\sum_{t=0}^{\infty}\frac{\Gamma\left(\frac{M}{2}+t\right)}{t!\left(\tr
\mathbf{\Sigma}^{-1}\mathbf{W}\mathbf{W}'\right)^{t}}\sum_{\kappa}
\frac{C_{\kappa}(R\mathbf{\Omega}
\mathbf{\Sigma}^{-1}\mathbf{W}\mathbf{W}')}{\left(\frac{1}{2}K\right)_{\kappa}}.
\end{eqnarray*}
So, we have proved that

\begin{cor}\label{th:QRreflectionshapeNORMAL}
The Gaussian QR reflection shape density is
\begin{eqnarray}\label{eq:QRreflectionshapeNORMAL}
f_{\mathbf{W}}(\mathbf{W})&=&\frac{\etr\{-\frac{1}{2}\mathbf{\Omega}\}\left(\tr
\mathbf{\Sigma}^{-1}\mathbf{W}\mathbf{W}'\right)^{-\frac{M}{2}}J(\mathbf{u})\displaystyle
\prod_{i=1}^{n}w_{ii}^{K-i}} {\pi^{\frac{M-nK}{2}}2^{-n+1}
\Gamma_{n}\left[\frac{1}{2}K\right]|\mathbf{\Sigma}|^{K/2}}\nonumber\\
&&\label{eq:QRreflectionshapeNORMAL}\times
\sum_{t=0}^{\infty}\frac{\Gamma\left(\frac{M}{2}+t\right)}{t!\left(\tr
\mathbf{\Sigma}^{-1}\mathbf{W}\mathbf{W}'\right)^{t}}\sum_{\kappa}
\frac{C_{\kappa}(\frac{1}{2}\mathbf{\Omega}
\mathbf{\Sigma}^{-1}\mathbf{W}\mathbf{W}')}{\left(\frac{1}{2}K\right)_{\kappa}},
\end{eqnarray}
where $M=(N-1)K$.
\end{cor}

The isotropic case of this density was derived by \citet{GM93}, and it is obtained
from (\ref{eq:QRreflectionshapeNORMAL}) noting that $C_{\kappa}(aB) =
a^{t}C_{\kappa}(B)$ and, if $\mathbf{\Sigma} = \sigma^{2}\mathbf{I}$, then
$\left(\tr\mathbf{\Sigma}^{-1}\mathbf{W}\mathbf{W}'\right)^{-\frac{M}{2}} =
\sigma^{M}$, $|\mathbf{\Sigma}|^{K/2} = \sigma^{M}$ and that
\begin{eqnarray*}
  \frac{C_{\kappa}(\frac{1}{2}\mathbf{\Omega} \mathbf{\Sigma}^{-1}\mathbf{W}\mathbf{W}')}{\left(\tr
   \mathbf{\Sigma}^{-1}\mathbf{W}\mathbf{W}'\right)^{t}} &=& \frac{C_{\kappa}\left(\displaystyle\frac{1}{2\sigma^{2}}
   \left(\frac{\boldsymbol{\mu}\boldsymbol{\mu}'}{\sigma^{2}}\right)\mathbf{W}\mathbf{W}'
   \right)}{\displaystyle\left(\frac{1}{\sigma^{2}}\right)^{t}}, \\
   &=& C_{\kappa}\left(\displaystyle\frac{1}{2\sigma^{2}}
        \boldsymbol{\mu}'\mathbf{W}\mathbf{W}' \boldsymbol{\mu} \right), \\
  &=& 2^{t}C_{\kappa}\left(\displaystyle\frac{1}{4\sigma^{2}}
        \boldsymbol{\mu}'\mathbf{W}\mathbf{W}' \boldsymbol{\mu} \right).
\end{eqnarray*}

Finally, we propose the result for the Kotz type I model
$$
  h(y)=\frac{R^{\tau-1+\frac{K(N-1)}{2}}\Gamma\left(\frac{K(N-1)}{2}\right)}{\pi^{K(N-1)/2}
  \Gamma\left(\tau-1+\frac{K(N-1)}{2}\right)}y^{\tau-1}\exp(-Ry),
$$

\begin{cor}\label{th:QRreflectionshapeKotz}
The Kotz type I QR reflection shape density is
\begin{eqnarray*}
f_{\mathbf{W}}(\mathbf{W})&=&\frac{\displaystyle
\prod_{i=1}^{n}w_{ii}^{K-i}J(\mathbf{u})}
{\Gamma_{n}\left[\frac{1}{2}K\right]|\mathbf{\Sigma}|^{K/2}}
\sum_{t=0}^{\infty}\sum_{\kappa}\frac{C_{\kappa}(R\mathbf{\Omega}
\mathbf{\Sigma}^{-1}\mathbf{W}\mathbf{W}')}{t!\left(\frac{1}{2}K\right)_{\kappa}}\nonumber\\
&&\times \frac{R^{\tau-1}\Gamma\left(\frac{M}{2}\right)\left(\tr
\mathbf{\Sigma}^{-1}\mathbf{W}\mathbf{W}'\right)^{\frac{M}{2}-t}\left(\tr
\mathbf{\Omega}\right)^{\tau-1}} {2^{-n+1}\pi^{\frac{M-nK}{2}}\etr
(R\mathbf{\Omega})}
\\&&\times\left\{\sum_{u=0}^{\infty}\frac{\Gamma\left(\frac{M}{2}+t+u\right)
\prod_{s=0}^{u-1}(\tau-1-s)}{u!R^{u}(\tr\mathbf{\Omega})^{u}
\Gamma\left[\tau-1+\frac{M}{2}\right] }\right.
\\&&\left.+
\sum_{m=1}^{k}\binom{k}{m}\left[\prod_{i=0}^{m-1}(\tau-1-i)\right]
\frac{(-R)^{-m}\left(\tr \mathbf{\Omega}\right)^{-m}}
{\Gamma\left[\tau-1-m+\frac{M}{2}\right]}\right.
\\&&\times\left.\sum_{u=0}^{\infty}\frac{\Gamma\left[\frac{M}{2}+t+u\right]
\prod_{s=0}^{u-1}(\tau-1-m-s)}{u!R^{u}(\tr\mathbf{\Omega})^{u} }\right\},
\end{eqnarray*}
where $M=(N-1)K$.
\end{cor}

\textit{Proof.} So, the corresponding $k$-th derivative follows from

$$
 \frac{d^{k}}{dy^{k}}[y^{T-1}\exp\{-Ry\}]=(-R)^{k}y^{\tau-1}\exp\{-Ry\}\left\{1+\sum_{m=1}^{k}\binom{k}{m}
 \left[\prod_{i=0}^{m-1}(\tau-1-i)\right](-Ry)^{-m}\right\}.\nonumber\\
$$
and the associating  QR reflection shape density can be obtained after some
simplification as

\begin{eqnarray*}
f_{\mathbf{W}}(\mathbf{W})&=&\frac{2^{n}\pi^{nK/2}\displaystyle
\prod_{i=1}^{n}w_{ii}^{K-i}J(\mathbf{u})}
{\Gamma_{n}\left[\frac{1}{2}K\right]|\mathbf{\Sigma}|^{K/2}}
\sum_{t=0}^{\infty}\sum_{\kappa}\frac{C_{\kappa}(\mathbf{\Omega}
\mathbf{\Sigma}^{-1}\mathbf{W}\mathbf{W}')}{t!\left(\frac{1}{2}K\right)_{\kappa}}\nonumber\\
&&\times\int_{0}^{\infty}r^{M+2t-1}h^{(2t)}[r^{2}\tr
\mathbf{\Sigma}^{-1}\mathbf{W}\mathbf{W}'+\tr\mathbf{\Omega}](dr)\\
&=&\frac{\displaystyle
\prod_{i=1}^{n}w_{ii}^{K-i}J(\mathbf{u})}{\Gamma_{n}\left[\frac{1}{2}K\right]|\mathbf{\Sigma}|^{K/2}}
\sum_{t=0}^{\infty}\sum_{\kappa}\frac{C_{\kappa}(R\mathbf{\Omega}
\mathbf{\Sigma}^{-1}\mathbf{W}\mathbf{W}')}{t!\left(\frac{1}{2}K\right)_{\kappa}}\nonumber\\&&\times
\frac{R^{\tau-1}\Gamma\left(\frac{M}{2}\right)\left(\tr
\mathbf{\Sigma}^{-1}\mathbf{W}\mathbf{W}'\right)^{\frac{M}{2}-t}\left(\tr
\mathbf{\Omega}\right)^{\tau-1}} {2^{-n+1}\pi^{\frac{M-nK}{2}}\etr
(R\mathbf{\Omega})}
\\&&\times\left\{\sum_{u=0}^{\infty}\frac{\Gamma\left[\frac{M}{2}+t+u\right]
\prod_{s=0}^{u-1}(\tau-1-s)}{u!R^{u}(\tr\mathbf{\Omega})^{u}
\Gamma\left[\tau-1+\frac{M}{2}\right] }\right.
\\&&\left.+ \sum_{m=1}^{k}\binom{k}{m}\left[\prod_{i=0}^{m-1}(\tau-1-i)\right]
\frac{(-R)^{-m}\left(\tr \mathbf{\Omega}\right)^{-m}}
{\Gamma\left[\tau-1-m+\frac{M}{2}\right]}\right.
\\&&\times\left.\sum_{u=0}^{\infty}\frac{\Gamma\left[\frac{M}{2}+t+u\right]
\prod_{s=0}^{u-1}(\tau-1-m-s)}{u!R^{u}(\tr\mathbf{\Omega})^{u} }\right\}. \qed
\end{eqnarray*}

\subsection{Example: Mouse Vertebra}\label{sub:mouse}

This classical application is studied in the gaussian case by \citet{DM98}. Here we
consider again the same model and contrasted it, via the modified BIC$^{*}$ criterion
( \citet{YY07}), with two non gaussian models.

Here we study three models, the gaussian shape (G), and the Kotz (K) model for
$\tau=2$ and $\tau=3$. The Gaussian shape density was given in
(\ref{eq:QRreflectionshapeNORMAL}) and the remaining shape distributions follows by
taking $\tau=2$ and $\tau=3$ in theorem \ref{th:QRreflectionshapeKotz}. However, we
need to simplify binomial series involved in the terms in braces, this can be done
straightforwardly but tedious by mathematical induction, the results are summarized
as follows.

Namely, the shape density associated to the Kotz model indexed by $\tau=2$,
$R=\frac{1}{2}$ (and $s=1$) is given by:
\begin{eqnarray*}
f_{\mathbf{W}}(\mathbf{W})&=&\frac{\displaystyle
\prod_{i=1}^{n}w_{ii}^{K-i}J(\mathbf{u})}{\pi^{\frac{M-nk}{2}}\Gamma_{n}\left(\frac{K}{2}\right)}
\etr\left(-\frac{\boldsymbol{\mu}'\boldsymbol{\mu}}{2\sigma^{2}}\right)\sum_{t=0}^{\infty}\frac{2^{n}}{M}\\
&&\times\left\{\left(\tr\frac{\boldsymbol{\mu}'\boldsymbol{\mu}}{2\sigma^{2}}-2t\right)
\Gamma\left[\frac{M}{2}+t\right]\right.
\\&&\left.+\Gamma\left[\frac{M}{2}+t+1\right]\right\}
\sum_{\kappa}\frac{C_{\kappa}\left(\frac{1}{2\sigma^{2}}\boldsymbol{\mu}'\mathbf{W}\mathbf{W}'
\boldsymbol{\mu}\right)}{t!\left(\frac{K}{2}\right)}.
\end{eqnarray*}
where $M=(N-1)K$.

And the corresponding density for the Kotz model $\tau=3$, is obtained as:

\begin{eqnarray*}
f_{\mathbf{W}}(\mathbf{W})&=&\frac{\displaystyle
\prod_{i=1}^{n}w_{ii}^{K-i}J(\mathbf{u})}{\pi^{\frac{M-nk}{2}}\Gamma_{n}\left(\frac{K}{2}\right)}
\etr\left(-\frac{\boldsymbol{\mu}'\boldsymbol{\mu}}{2\sigma^{2}}\right)
\sum_{t=0}^{\infty}\frac{2^{n+1}}{M(M+2)}\\
&&\times\left\{\left(\tr\frac{\boldsymbol{\mu}'\boldsymbol{\mu}}
{2\sigma^{2}}-2t\right)\Gamma\left[\frac{M}{2}+t\right]+
2\left(\tr\frac{\boldsymbol{\mu}'\boldsymbol{\mu}}{2\sigma^{2}}-2t\right)
\Gamma\left[\frac{M}{2}+t+1\right]\right.
\\&&\left.+\Gamma\left[\frac{M}{2}+t+2\right]\right\}
\sum_{\kappa}\frac{C_{\kappa}\left(\frac{1}{2\sigma^{2}}\boldsymbol{\mu}'
\mathbf{W}\mathbf{W}'\boldsymbol{\mu}\right)}{t!\left(\frac{K}{2}\right)}.
\end{eqnarray*}
where $M=(N-1)K$.

In order to decide which the elliptical model is the best one, different criteria
have been employed for the model selection. We shall consider a modification of the
BIC$^{*}$ statistic as discussed in \citet{YY07}, and which was first achieved by
\citet{ri:78} in a coding theory framework. The modified BIC$^{*}$ is given by:
$$
  BIC^{*}=-2\mathfrak{L}(\widetilde{\boldsymbol{\mu}},\widetilde{\sigma}^{2},h)
         + n_{p}(\log(n+2 )- \log 24),
$$
where $\mathfrak{L}(\widetilde{\boldsymbol{\mu}},\widetilde{\sigma}^{2},h)$ is the
maximum of the log-likelihood function, $n$ is the sample size and $n_{p}$ is the
number of parameters to be estimated for each particular shape density.

As proposed by \citet{kr:95} and \citet{r:95}, the following selection criteria have
been employed for the model selection.

\medskip

\begin{table}[ht]  \centering \caption{Grades of evidence
corresponding to values of the $BIC^{*}$ difference.}\label{table2}
\medskip
\renewcommand{\arraystretch}{1}
\begin{center}
  \begin{tabular}{cl}
    \hline
    $BIC^{*}$ difference & Evidence\\
    \hline
    0--2 & Weak\\
    2--6 & Positive \\
    6--10 & Strong\\
    $>$ 10 & Very strong\\
    \hline
  \end{tabular}
\end{center}
\end{table}

The maximum likelihood estimators  for location and scale parameters associated with
the small and large groups are summarized in the following table:

\medskip

\begin{table}[ht!]
  \centering
  \caption{The maximum likelihood estimators}\label{tb:1}
  \medskip
\begin{center}
\begin{small}\label{Tab4MR}
\begin{tabular}{|c|c|c|c|c|c|c|c|}
  \hline
  Group& $BIC^{*}$ & $\widetilde{\mu}_{11}$ &$\widetilde{\mu}_{12}$& $\widetilde{\mu}_{21}$ &
    $\widetilde{\mu}_{22}$ & $\widetilde{\mu}_{31}$ & $\widetilde{\mu}_{32}$  \\
   & $\build{K:\tau=2}{G}{K:\tau=3}$&&&&&&\\
  \hline
  Small & $\build{-418.011}{-403.824}{-307.863}$ & $\build{-3.3846}{1.2398}{3.5716 }$ &
  $\build{44.7126 }{39.2181}{131.3120 }$
   & $\build{14.7682}{13.3663}{ 44.6939 }$
  &$\build{ 5.5268}{3.4263}{ 11.6686 }$    &$\build{ 25.3360}{22.1414}{ 74.1405 }$   &
  $\build{ 1.0451 }{-1.4618}{ -4.5674 }$
  \\ \hline
  Large& $\build{206.7321}{199.6375}{151.6613}$ &$\build{-7.2450 }{16.9915}{-26.5182}$
&$\build{ -90.9671 }{-104.1137}{-71.4992}$ &$\build{ 28.0714
}{34.6059}{20.0230}$
  & $\build{ -11.2058 }{-4.8256}{-15.3962}$   & $\build{ 58.7553 }{65.7152}{47.4424}$ &
  $\build{ 0.8674 }{17.2009}{-12.6667}$  \\
  \hline
\end{tabular}
\end{small}
\end{center}
\end{table}

\begin{center}
\begin{footnotesize}
\begin{tabular}{|c|c|c|c|c|c|c|c|}
  \hline
   $\widetilde{\mu}_{41}$ & $\widetilde{\mu}_{42}$ &
   $\widetilde{\mu}_{51}$ & $\widetilde{\mu}_{52}$
   & $\widetilde{\sigma^{2}}$  \\
  \hline
   $\build{5.2270}{4.0894}{13.7605 }$ &$\build{ -4.8965}{-4.7493}{-15.8399}$   &
   $\build{ -30.7140}{-27.1075}{ -90.7380}$
   &$\build{ -4.1166}{-0.7072}{  -2.7664 }$
     & $\build{ 48.6680}{42.8290}{289.4148}$
  \\ \hline
  $\build{ 7.1216 }{5.1349}{8.0814}$ &$\build{10.3519 }{13.5872}{6.7039}$   &
  $\build{ -72.3213 }{-82.9587}{-56.6913}$
   &$\build{6.4313 }{-12.7549}{21.6223}$ &
  $\build{225.3525}{346.0959}{109.0523}$\\
  \hline
\end{tabular}
\end{footnotesize}
\end{center}

According to the modified BIC$^{*}$ criterion (see Table \ref{Tab4MR}), the Kotz
model with parameters $\tau=2$, $R=\frac{1}{2}$ and $s=1$ is the most appropriate
among the three elliptical densities for modeling the data. There is a very strong
difference between the non gaussian and  the classical gaussian model widely detailed
and applied by \citet{DM98} (and previous works) in this experiment.

Let $\boldsymbol{\mu}_{1}$ and $\boldsymbol{\mu}_{2}$ be the mea shape of the small
and large groups, respectively. We test equal mean shape under the best model, and
the likelihood ratio (based on $-2\log\Lambda\approx\chi_{10}^{2}$) for the test
$H_{0}:\boldsymbol{\mu}_{1}=\boldsymbol{\mu}_{2}$ vs
$H_{a}:\boldsymbol{\mu}_{1}\neq\boldsymbol{\mu}_{2}$, provides the p-value $0.3\,
10^{-12}$, which means that there extremely evidence that the mean shapes of the two
groups are different.  This is the same conclusion obtained by \citet{DM98} for a
gaussian model.

A final comment, for any elliptical model we can obtain the Q reflection model,
however a nontrivial problem appears, the $2t$-th derivative of the generator model,
which can be seen as a partition theory problem. For the general case of a Kotz model
($s \neq 1$), and another models like Pearson II and VII, Bessel, Jensen-logistic, we
can use formulae for these derivatives given by \citet{Caro2009}. The resulting
densities have again a form of a generalized series of zonal polynomials which can be
computed efficiently after some modification of existing works for hypergeometric
series (see \citet{KE06}), thus the inference over an exact density can be performed,
avoiding the use of any asymptotic distribution, and the initial transformation
avoids the invariant polynomials of \citet{d:80}, which at present seems can not be
computable.

\section*{Acknowledgment}

This research work was supported  by University of Medellin (Medellin, Colombia) and
Universidad Aut\'onoma Agraria Antonio Narro (M\'{e}xico),  joint grant No. 469,
SUMMA group. Also, the first author was partially supported  by CONACYT - M\'exico,
research grant no. \ 138713 and IDI-Spain, Grants No. \ FQM2006-2271 and
MTM2008-05785 and the paper was written during J. A. D\'{\i}az- Garc\'{\i}a's stay as
a visiting professor at the Department of Statistics and O. R. of the University of
Granada, Spain.


\begin{thebibliography}{}

\bibitem[Caro-Lopera \textit{et al.} (2009)]{Caro2009}
    F. J. Caro-Lopera, J. A. D\'{\i}az-Garc\'{\i}a and G. Gonz\'{a}lez-Far\'{\i}as,
    Noncentral elliptical configuration density,
    J. Multivariate Anal. 101(1) (2009), 32--43.

\bibitem[Davis(1908)]{d:80}
    A. W. Davis,
    Invariant polynomials with two matrix arguments, extending the zonal polynomials,
    in: \emph{Multivariate Analysis V}, (Krishnaiah, P. R. ed.), North-Holland, 1980.

\bibitem[D\'{\i}az-Garc\'{\i}a and Gonz\'{a}lez-Far\'{\i}as (2005)]{dggg05}
    J. A. D\'{\i}az-Garc\'{\i}a and G. Gonz\'{a}lez-Far\'{\i}as,
    Singular random matrix decompositions: Distributions,
    J. Multivariate Anal. 194(1) (2005), 109--122.

\bibitem[D\'{\i}az-Garc\'{\i}a and Guti\'errez-J\'aimez (2006)]{dggj06}
    J. A. D\'{\i}az-Garc\'{\i}a and R. Guti\'errez-J\'aimez,
    Wishart and Pseudo-Wishart distributions
    under elliptical laws and related distributions in the shape theory context,
    J. Stat. Plan. Inference 136(12) (2006), 4176--4193.

\bibitem[Dryden and Mardia (1998)]{DM98}
    I. L. Dryden and K.V. Mardia,
    Statistical shape analysis,
    John Wiley and Sons, Chichester, 1998.

\bibitem[Fang and Zhang (1990)]{fz:90}
    K. T. Fang,  and Y. T. Zhang, Generalized
    Multivariate Analysis,
    Science Press, Springer-Verlag, Beijing, 1990.

\bibitem[Goodall and Mardia(1991)]{gm:91}
   C. R. Goodall, and K. V. Mardia,
   A geometrical derivation of the shape density,
   Adv. in Appl. Probab. 23 (1991) 496--514.

\bibitem[Goodall and Mardia (1993)]{GM93}
    C. R. Goodall, and K. V. Mardia,
    Multivariate Aspects of Shape Theory,
    Ann.  Statist. 21 (1993) 848--866.

\bibitem[Gupta and Varga(1993)]{gv:93}
    A. K. Gupta, and T. Varga,
    Elliptically Contoured Models in Statistics,
    Kluwer Academic Publishers, Dordrecht, 1993.

\bibitem[James(1964)]{JAT64}
    A. T. James,  Distributions of matrix variate and latent roots
   derived from normal samples,
   Ann. Math. Statist. 35 (1964) 475--501.

\bibitem[Kass and Raftery (1995)]{kr:95}
    R. E. Kass, and A. E. Raftery, Bayes factor, J. Amer. Statist. Soc. 90 (1995) 773--795.

\bibitem[Koev and Edelman (2006)]{KE06}
    P. Koev and A. Edelman,
    The efficient evaluation of the hypergeometric function of a matrix argument,
    Math. Comp. 75 (2006) 833--846.

\bibitem[Muirhead(1982)]{MR1982}
    R. J. Muirhead,
    Aspects of multivariate statistical theory,
    Wiley Series in Probability and Mathematical Statistics. John Wiley \& Sons, Inc.
    1982.

\bibitem[Raftery(1995)]{r:95}
    A. E. Raftery, Bayesian model selection in social research,
    Sociological Methodology, 25 (1995) 111--163.

\bibitem[Rissanen(1978)]{ri:78}
    J. Rissanen, Modelling by shortest data description,
    Automatica, 14 (1978) 465--471.

\bibitem[Yang and Yang(2007)]{YY07}
    Ch. Ch. Yang and Ch. Ch. Yang,
    Separating latent classes by information criteria,
    J. Classification 24 (2007) 183--203.

\end{thebibliography}
\end{document}